\documentclass[11pt]{article}
\usepackage{graphicx}
\usepackage[dvipsnames]{color}
\definecolor{gris}{gray}{0.7}
\usepackage{color}
\definecolor{nar}{rgb}{.8,0,1}
\definecolor{otro}{rgb}{0,1,1}
\usepackage{amssymb}
\usepackage{amsmath}
\usepackage{epstopdf}
\usepackage{multicol}
\usepackage{epic}
\usepackage{fancybox}
\usepackage{url}
\usepackage{bbm}
\usepackage{dsfont}
\usepackage{eurosym}
\usepackage[applemac]{inputenc}
\usepackage[dvipsnames]{color}
\definecolor{gris}{gray}{0.7}

\newcommand{\D}{\displaystyle}
\newcommand{\R}{\mathbb{R}}
\newcommand{\N}{\mathbb{N}}

\newcommand{\C}{\mathbb{C}}

\newtheorem{theorem}{Theorem}
\newtheorem{remark}{Remark}

\newtheorem{claim}{Claim}
 \numberwithin{equation}{section}

\begin{document}
  
\title{Lower Bounds for Non-Trivial Traveling Wave Solutions of Equations of KdV Type}

\author{C.E. Kenig, G. Ponce, L. Vega\footnote{
The first and second authors  are supported by NSF grants DMS--0968742 and DMS--1101499 respectively.
The third author by MEC grant MTM 2007--03029. }
\\\\
C.E. Kenig: Department of Mathematics, University of Chicago, \\
Chicago, IL 60637, USA.\\
cek@math.uchicago.edu\\\\
G. Ponce: Department of Mathematics, University of California,\\
Santa Barbara, CA 93106, USA.\\
ponce@math.ucsb.edu\\\\\
L. Vega\\
L. Vega: Departamento de Matem\'aticas, Universidad del Pa\'is Vasco/EHU, \\
 Apartado 644, 48080, Bilbao, Spain\\
luis.vega@ehu.es
}

\maketitle

\begin{abstract} We prove that if a solution of an equation of KdV type is bounded above by a traveling wave with an amplitude that decays faster than a given linear exponential then it must be zero. We assume no restrictions neither on the size nor in the direction of the speed of the traveling wave.
\end{abstract}

 \section{Introduction and main results.}
  \label{S1}

In this paper we will continue our study initiated in \cite{EKPV1} on solutions of dispersive systems that decay like a linear exponential. In particular, we extend the results proved in that paper for traveling wave solutions of non--linear Schr\"odinger (NLS) equations to the case of equations of Korteweg--de Vries (KdV) type.

More concretely we shall consider solutions of equations of the form
\begin{equation}
\label{eq1.1} \left(\partial_t+\partial_x^3\right)u=a(u)\partial_x u,\qquad x\in\R\,\,,\,\, t\in\R,
\end{equation}
with $a$ regular and such that
\begin{equation}
\label{eq1.2} \left|a(s)\right|\le M_1\left(|s|+|s|^j\right),\qquad j=1,2,3,\dots
\end{equation}

Our main result is the following one.

 \begin{theorem}
  \label{Th1}
  Assume $u\in\mathcal C^1\left(\R\,:\,H^1(\R)\right)$ is a real solution of \eqref{eq1.1}--\eqref{eq1.2} with $a(u)$ also real. Then there exist $\lambda_0\ge1$ and $c_0>0$ such that if for some $\lambda>0$, $b\in\R$
 \begin{eqnarray}
\label{eq1.3} \sup_{t\ge 0}\left\|u(\cdot,t)\right\|_{H^1}^2\le M_2,\\
\nonumber\\
\label{eq1.4}\sup_{t\ge 0}\int e^{|x-bt|}\left|u(x,t)\right|^2dx\le M_3,\\
\nonumber\\
\label{eq1.5}\sup_{t\ge 0}\int e^{2\lambda|x-bt|}\left|u(x,t)\right|^2<+\infty,
\end{eqnarray} 
and 
$$\lambda\ge \max\left\{\lambda_0,c_0M_1^2\left(1+\left[M_3+M_2^{1/2}M_3^{1/2}\right]^{j/2}\right)\right\},$$
then
$$u\equiv 0.$$
 \end{theorem}
 
 Relevant examples of non--linear potentials are the pure power ones
 \begin{equation}
 \label{eq1.6} a(u)=-ju^j\qquad j=1,2,3...,
 \end{equation}
 and the completely integrable ones
 \begin{equation}
\label{eq1.7}
a(u)=\alpha_1u+\alpha_2u^2,
\end{equation}
that are usually known as Gardner equations, where the cases $(\alpha_1, \alpha_2)=(1,0)$ and $(\alpha_1, \alpha_2)=(0,1)$ are the KdV  and modified KdV equation respectively.

 \begin{remark}
  \label{remark1} For  \eqref{eq1.6} there is a scaling symmetry 
$$u_{\lambda}(x,t)=\lambda^{2/j} u\left(\lambda x,\lambda^3t\right),\quad\lambda>0,$$
and  an explicit family of traveling waves
$$\lambda^{2/j}\varphi_j\left(\lambda(x-\lambda^2t)\right),\qquad\text{with }\, \varphi_j=\left(\frac{j+2}{2}\text{sech}^2(\frac{j}{2}x)\right)^{\frac{1}{j}}.$$
 Hence in this case we have traveling wave solutions that propagate to the right with velocity $\lambda^2$. 
 
Also given any $j$ and using the scaling invariance  we can always choose $\lambda$  small enough such that   there exists a non-trivial solution that verifies \eqref{eq1.3}, \eqref{eq1.4}, and \eqref{eq1.5}.
\end{remark}
Regarding negative velocities we recall the existence of breather solutions in the case of Gardner equations for some specific choices of the parameter $\alpha_1$ and $\alpha_2$. Let us consider for simplicity the modified KdV equation; and specifically we take $a(u)=-6u^2$. Then, with this normalization the soliton solutions are\begin{equation}
\label{eq1.8}
u(x,t)=\lambda\,\text{sech}\,\big(\lambda(x-\lambda^2t)\big),
\end{equation}
and the breather solutions are \cite{Wa}
\begin{equation}
\label{eq1.9}
u(x,t)=2\lambda\,\text{sech}\,\big(\lambda(x+\gamma t)\big)\cdot\left\{\D\frac{\cos\left[\Phi(x,t)\right]-\frac{\lambda}{\mu}\sin \Phi(x,t)\tanh\big(\lambda(x+\gamma t)\big)}
{1+\left(\frac{\lambda}{\mu}\right)^2\sin^2\left(\Phi(x,t)\right)\text{sech}^2\,\big(\lambda(x+\gamma t)\big)}\right\},
\end{equation}
with
$$\gamma=3\mu^2-\lambda^2\quad,\quad
\delta=\mu^2-3\lambda^2\quad,\text{and}\quad
\Phi=\mu(x+\delta t)-\tan^{-1}\left(\D\frac{\lambda}{\mu}\right).$$

These breathers can be seen as wave packets with an amplitude and a frequency determined by  the parameters $\lambda$
 and $\mu$. The velocity of the amplitude is $\gamma=3\mu^2-\lambda^2$ and it can take any real value.
 
 \begin{remark}
  \label{remark2}
  Notice that in the above examples the decay rate of the amplitude is just given by $\lambda$ as in the statement of our theorem. Choosing as before $\lambda$ small enough  we also conclude the sharpness of our result for negative velocities.
\end{remark}

\begin{remark}
  \label{remark3}
Breather solutions also exist for complex KdV equations. For example if we choose $a(u)=-|u|^2$ with $u(x,t)\in\C$ it is easy to check (see for example \cite{KPV2}) that
$$u(x,t)=\sqrt 3e^{-it\mu(3\lambda^2-\mu^2)+i\mu x} \psi_{\lambda}\left(x-(\lambda^2-3\mu^2)t\right),$$
with $\psi_{\lambda}(x)=\lambda\psi_1(\lambda x)$ and $\psi_1(x)=\sqrt 2 \text{sech}\,x$. The arguments needed for the proof of Theorem \ref{Th1} also work for this complex KdV (i.e. $a(u)=-|u|^2$).
\end{remark}

There are three fundamental ingredients in the proof of Theorem \ref{Th1}. The first one is that the $L^2$ norm is a conserved quantity. This property is needed to conclude that $u(x,0)=0$ if for some sequence of times $t_n$ we have that
$$\lim_{t_n\to\infty}\|u(\cdot,t_n)\|_{L^2}=0.$$

The second ingredient is Kato's theory in \cite{Ka} concerning the persistent property of solutions to  the initial value problem with data that satisfies linear exponential decay. Kato's result is used to prove that if \eqref{eq1.5} holds then a similar inequality is also true for $\partial_x^k u\,$ $\, k=1,2,3$. 
We require this property since instead of working with the equation \eqref{eq1.1}--\eqref{eq1.2} we use the one satisfied by $f=e^{\lambda\theta}u$ with $\theta(x,t)$ an appropriate Carleman weight that grows almost linearly at infinity. 

The final ingredient is  the convexity of $H(t)=\|f\|_{L^2}^2=\langle f,f\rangle$. As done in \cite{EKPV1} and \cite{EKPV2}  we estimate
$\dot{H}(t)$ and $\ddot{H}(t)$ using integration by parts, therefore we need $f$ and its spatial derivatives in $L^2$. At this point we closely follow the arguments in \cite{EKPV1} but the algebra, as it can be  expected, turns out to be more complicated in the KdV setting than in the NLS one. The details can be found in the proof of Claim \ref{CL1} in Section \ref{S2}. In this claim we establish the positivity of the commutator that appears in the computation of  $\ddot{H}(t)$ -see also \cite{R} and \cite{SS}. 

There is still another difference with respect to \cite{EKPV1} . For generalized KdV \eqref{eq1.1}-\eqref{eq1.2} the non--linear potential can not be treated in a perturbation way and some structure is needed (see the proof of Claim \ref{CL2} in Section  \ref{S2}). This structure still holds in the complex KdV equation that we mentioned above.

The use of Carleman weights to obtain positive commutators is  a standard technique in elliptic theory for obtaining lower bounds of eigenfunctions of a Schr\"odinger operator  -see \cite{C}, \cite{FH}, \cite{M}. In the time evolution setting some modifications of this technique are needed and an important step in our approach is the use of the two identities \eqref{eq2.4} and \eqref{eq2.5}. These  identities  also appear in \cite{EKPV1}.

The rest of the paper, that is to say Section  \ref{S2}, is devoted to the proof of Theorem \ref{Th1}.

 \section{Proof of Theorem \ref{Th1}.}
  \label{S2}
  
  From the results by T. Kato in \cite{Ka} we can assume that \eqref{eq1.5} is also satisfied possibly with another $\lambda$ for $\partial_x^k u(x,t)$ with $k=1,2,3$.
  
  As in \cite{EKPV1} we shall work with $f=e^{\lambda\theta}u$ with $\theta=\theta(x,t)$. More concretely $\theta(x,t)=\varphi(r)$,  $r=|x-bt|$ for some regular even $\varphi$ that will be fixed later on and that it grows at most linearly at infinity. As in \cite{EKPV2} we have
  \begin{eqnarray}
\nonumber &e ^{\lambda\theta}\left(\partial_t+\partial_x^3\right)e^{-\lambda\theta}f=(S_{\lambda}+A_{\lambda})f;\\
\nonumber &S_{\lambda}=-3\lambda\partial_x\left(\partial_x\theta\partial_x\right)+\left(-\lambda^3(\partial_x\theta)^3-\lambda\partial_x^3\theta-\lambda\partial_t\theta\right);\\
\nonumber &A_{\lambda}=\partial_t+\partial_x^3+3\lambda^2\left(\partial_x\theta\right)^2\partial_x+3\lambda^2\partial_x\theta\partial_x^2\theta= \partial_t+\tilde A_{\lambda};\\
\label{eq2.0} &A_{\lambda}^{\ast}=-A_{\lambda}\qquad \tilde A_{\lambda}^{\ast}=-\tilde A_{\lambda}\qquad S_{\lambda}^{\ast}=S_{\lambda}.
\end{eqnarray}

Hence, using the notation $\partial_x f=f_x,$

\begin{equation}
\label{eq2.1}
\begin{array}{rcl}
\langle[S_{\lambda};A_{\lambda}]f,f\rangle&=&\langle (S_{\lambda} A_{\lambda}-A_{\lambda}S_{\lambda})f,f\rangle\\\\
&=&9\lambda\D\int\varphi^{\prime\prime}(r)f_{xx}^2+
6\lambda b\int \varphi^{\prime\prime}(r)f_x^2+
18\lambda^3\int\left(\varphi^{\prime}(r)\right)^2\varphi^{\prime\prime}(r)f_x^2\\\\
&-&
6\lambda\int\varphi^{IV}(r)f_x^2+9\lambda^5\D\int\left(\varphi^{\prime}(r)\right)^4\varphi^{\prime\prime}(r)f^2+
\lambda b^2\int \varphi^{\prime\prime}(r)f^2\\\\
&-&3\lambda^3\D\int\left(\varphi^{\prime\prime}(r)\right)^3f^2-
18\lambda^3\int \varphi^{\prime}(r)\varphi^{\prime\prime}(r)\varphi^{\prime\prime\prime}(r)f^2\\\\
&-&3\lambda^3\D\int\left(\varphi^{\prime\prime}(r)\right)^2\varphi^{IV}(r)f^2+
\lambda\int\varphi^{VI}(r)f^2-
2\lambda b\int\varphi^{IV}(r)f^2\\\\
&-&6\lambda^3 b\D\int\varphi^{\prime}(r)^2\varphi^{\prime\prime}(r)f^2\\\\
&=&\textcircled{\tiny{1}}+\textcircled{\tiny{2}}+\textcircled{\tiny{3}}+
\textcircled{\tiny{4}}+\textcircled{\tiny{5}}+\textcircled{\tiny{6}}+
\textcircled{\tiny{7}}+\textcircled{\tiny{8}}+\textcircled{\tiny{9}}+
\textcircled{\tiny{10}}+\textcircled{\tiny{11}}+\textcircled{\tiny{12}}.
\end{array}
\end{equation}

Notice that
\begin{equation}
\label{eq2.2}
\textcircled{\tiny{5}}+\textcircled{\tiny{6}}+\textcircled{\tiny{12}}=\lambda\D\int\left(3\lambda^2\left(\varphi^{\prime}(r)\right)^2-b\right)^2\varphi^{\prime\prime}(r)f^2.
\end{equation}

As in \cite{EKPV1} we shall use two identities that hold for solutions of
\begin{equation}
\label{eq2.3}
\partial_t f=-\left(S_{\lambda}+\tilde A_{\lambda}\right)f+F,
\end{equation}
with
$$\tilde A_{\lambda}^{\ast}=-\tilde A_{\lambda}\qquad S_{\lambda}^{\ast}=S_{\lambda}.$$
First observe  that from \eqref{eq2.3}
$$\frac{d}{dt} \langle f,f\rangle=-2 \langle S_{\lambda}f,f\rangle +2\langle F,f\rangle,$$
and, recalling that in our case $A_{\lambda}=\partial_t+\tilde A_{\lambda}$, we have
$$\frac{d}{dt} \langle S_{\lambda}f,f\rangle= -\langle (S_{\lambda}A_{\lambda}-A_{\lambda}S_{\lambda})f,f\rangle 
-2 \langle S_{\lambda}f,S_{\lambda}f\rangle+2\langle F,S_{\lambda}f\rangle.$$
In this last identity we have used that
$$S_{\lambda}A_{\lambda}-A_{\lambda}S_{\lambda}=S_{\lambda}\tilde A_{\lambda}-\tilde A_{\lambda}S_{\lambda}-(S_{\lambda})_t.$$

Then, take $\eta:[a,b]\longrightarrow \R$. A simple integration by parts   gives (see Proposition 1 and Proposition 2 in \cite {EKPV1} for more details)
\begin{equation}
\label{eq2.4}
\D\int_a^b\eta^{\prime}(t)\langle S_{\lambda}f,f\rangle \,dt=
-\frac 12\left.\left(\eta^{\prime}\langle f,f\rangle\right)\right|_a^b
+\int_a^b\eta^{\prime}(t)\langle F,f\rangle dt
+\frac 12\int_a^b\eta^{\prime\prime}(t)\langle f,f\rangle dt,
\end{equation}
and
\begin{equation}
\label{eq2.5}
\begin{array}{rcl}
\D\int_a^b\eta^{\prime}(t)\langle S_{\lambda}f,f\rangle \,dt&=&
\left.\left(\eta\langle S_{\lambda}f,f\rangle\right)\right|_a^b
+\int_a^b\eta\langle (S_{\lambda}A_{\lambda}-A_{\lambda}S_{\lambda}f),f\rangle dt
\\\\
&+&
2\int_a^b\eta\langle S_{\lambda}f,S_{\lambda} f\rangle dt
-2\int_a^b\eta\langle F, S_{\lambda}f\rangle dt.
\end{array}
\end{equation}

Notice that if $u$ solves \eqref{eq1.1}--\eqref{eq1.2} and $f=e^{\lambda\varphi(x-bt)}u$, then $f$ solves
\begin{equation}
\label{eq2.6}
\partial_t f=
-\left(S_{\lambda}+\tilde A_{\lambda}\right)f+a(u)\left(\partial_x f-\lambda\varphi^{\prime}(x-bt)f\right).
\end{equation}

Finally we shall use the following Carleman weight $\varphi_0$:
\begin{equation}
\label{eq2.7}
\varphi_0\in\mathcal C^6(\R), \text{ even and positive,}
\end{equation}
\begin{equation}
\label{eq2.8}
\varphi^{\prime}_0(r)=r\,\text{ if }\, 0\le r\le 3/2\,\text{ and }\,\varphi_0^{\prime}(r)=2-\D\frac {\log 2}{4\log r}\,\text{ if }\,r\ge 2,
\end{equation}
\begin{equation}
\label{2.9}
0<\varphi_0^{\prime\prime}(r)\le 1\,\text{ and it is a decreasing function for }\, r>3/2,
\end{equation}
\begin{eqnarray}
\nonumber&\text{there exists $c_0>0$ such that}\\
\nonumber\\
\label{eq2.10}&
\left|\D\frac{d^k}{dr^k}\varphi_0(r)\right|\le c_0\varphi_0^{\prime\prime}(r)\qquad k=3,4,5,6.
\end{eqnarray}

The proof of the theorem will follow from the next three claims. Recall that $f=e^{\lambda\varphi(r)}u$ with $r=|x-bt|$. If $b\geq\frac{3}{2}\lambda^2$ we shall choose $\varphi=\frac{1}{4}\varphi_0$. Otherwise we shall take $\varphi=\varphi_0$. The reason behind these choices will become clear in the proof of Claim \ref{CL1} where different cases will be considered.

\begin{claim}
\label{CL1}
There exist $\lambda_0\ge 1$, $A_0>0$ such that for all $\lambda\ge\lambda_0$
\begin{equation}
\label{eq2.11}
\langle\left(S_{\lambda}A_{\lambda}-A_{\lambda}S_{\lambda}\right)f,f\rangle\ge
A_0\D\int
\varphi^{\prime\prime}_0\left(\lambda^3f^2+\lambda^2f^2_x\right).
\end{equation}
\end{claim}

The proof of this claim is long and it is postponed.

\begin{claim}
\label{CL2}
There exists a constant $C>0$ which depends on \eqref{eq1.5} and  on $M_1$, $M_2$ and $M_3$ given in \eqref{eq1.2}, \eqref{eq1.3}, and \eqref{eq1.3} such that
\begin{equation}
\label{eq2.12}
\left|\langle a(u)\left(\partial_x f-\lambda\varphi^{\prime}(x-bt)f\right),f\rangle\right|\le \lambda C.
\end{equation}
\end{claim}

 \begin{remark}
  \label{remark3} In the above inequality \eqref{eq2.12} the structure of the non-linear term plays a role. This makes a difference with respect to the non-linear Schr\"odinger equation, see \cite{EKPV1}.
  \end{remark}

\underline{Proof} Recall that
$$u=e^{-\lambda\varphi(x-bt)}f,$$
so that from \eqref{eq1.3} $u$ is in $L^2$. Also observe that
$$e^{\lambda\varphi}\partial_x u=\partial_x f-\lambda \varphi^{\prime} f.$$
Hence by integration by parts
$$\D\int a(u)\left(\partial_xf-\lambda\varphi^{\prime}f\right)f\,dx=
\int a(u)e^{\lambda\varphi}\partial_x u \,e^{\lambda \varphi}u\,dx=
-\int e^{2\lambda\varphi}2\lambda\varphi'u^2\alpha(u)\,dx,$$
with
$$\alpha(s)=\D\frac 1{s^2}\int_0^s a(s^{\prime})s^{\prime}ds^{\prime}.$$
The claim easily follows from the boundedness of   $\varphi_0^{\prime}$ and \eqref{eq1.2}-\eqref{eq1.5}.

\begin{claim}
\label{CL3}
There exists a universal constant $C_1$ such that
\begin{equation}
\label{eq2.13}
\left\|a(u)\left(f_x-\lambda\varphi^{\prime}f\right)\right\|_{L^2}^2\le
C_1\D\frac{M_1^2}{\lambda}
\left[\left(M_2^{1/2}M_3^{1/2}+M_3\right)^{j/2}+1\right]
\int\varphi^{\prime\prime}_0\left(\lambda^3f^2+\lambda f_x^2\right)dx.
\end{equation}
\end{claim}

\underline{Proof}  We have
\begin{equation}
\label{eq2.133}
\begin{array}{rcl}
e^{\frac 12|x-bt|}u^2&\le&
\D\int_{-\infty}^{\infty}\left|\frac d{dy}\left( e^{\frac 12|y-bt|}u^2(y)\right)\right|dy\\\\
&\le&2\|u_x\|_{L^2}\left\|e^{\frac 12|y-bt|}u\right\|_{L^2}+
\D\frac 12\left\|e^{\frac 12|y-bt|}u\right\|_{L^2}^2\\\\
&\le& 2\left(M_2M_3\right)^{1/2}+\D\frac 12M_3.
\end{array}
\end{equation}

Also there is a universal constant $C_0$ such that
\begin{equation}
\label{eq2.1333}e^{-\frac 12 r}\le C_0\varphi^{\prime\prime}_0(r).
\end{equation}
Then we get from \eqref{eq2.133}
\begin{equation}
\label{eq2.1334}|a(u)|^2\leq C_0 \varphi_0^{\prime\prime}( 2\left(M_2M_3\right)^{1/2}+\D\frac 12M_3).
\end{equation}

As we already said, either $\varphi=\frac{1}{4}\varphi_0$, or $\varphi=\varphi_0$. Hence
\begin{equation}
\label{eq2.13344}|f_x-\lambda\varphi^{\prime}f|^2\leq 2|f_x|^2+2\lambda^2\sup\varphi_0^{\prime}|f|^2\leq C_0(|f_x|^2+\lambda^2|f|^2),
\end{equation}
for another universal constant $C_0$. Also recall that $\lambda\geq1$.
The claim follows from  \eqref{eq2.1334} and \eqref{eq2.13344}.

Let us finish the proof of the theorem before we prove Claim \ref{CL1}. We follow the arguments in \cite{EKPV1}, and more concretely those in Proposition 1 and Proposition 2 of that paper.

 Take $\lambda\geq\lambda_0$ and such that
$$C_1\D\frac{M_1^2}{\lambda}
\left[\left(M_2^{1/2}M_3^{1/2}+M_3\right)^{j/2}+1\right]\le\frac{A_0}4.$$
Recall \eqref{eq2.6} so that \eqref{eq2.3} is satisfied with
$$F= a(u) (f_x-\lambda\varphi^{\prime}(x-bt)f).$$
Then, from \eqref{eq2.11} and \eqref{eq2.13} we get
\begin{equation}
\label{eq2.14}
\langle\left[S_{\lambda};A_{\lambda}\right]f,f\rangle+
2\langle S_{\lambda}f,S_{\lambda}f\rangle-
2\langle F, S_{\lambda}f\rangle\ge
\langle  S_{\lambda}f,S_{\lambda}f\rangle+\frac{A_0}{2}\D\int
\varphi^{\prime\prime}_0\left(\lambda^3f^2+\lambda^2f^2_x\right).
\end{equation}

Our first step is to use the above inequality to obtain a uniform estimate for $\langle  S_{\lambda}f,S_{\lambda}f\rangle$.  Take for any $n\in\N$,  $\eta_n(t)$ defined on $\left[n-1/2,n+1/2\right]$ as 
$$\eta_n(t)=\D\frac 12-|t-n|.$$
Notice that $\eta_n(n\pm1/2)=0$, so that the right hand side of \eqref{eq2.5} is lower bounded using  \eqref{eq2.14} by
$$\D\int\eta_n\left\|S_{\lambda}f\right\|_{L^2}^2dt.$$
On the other hand $|\eta^{\prime}_n|=1$
and  $\eta_n^{\prime\prime}=-\delta(t-n)\leq 0$. Then, using \eqref{eq1.5} and \eqref{eq2.12}, the right hand side of \eqref{eq2.4} is bounded above by $\lambda C$ with $C$ depending on $M_1,M_2, M_3$, and the bound in  \eqref{eq1.5}, but not on $n$.

As a consequence we get
$$\D\int\eta_n\left\|S_{\lambda}f\right\|_{L^2}^2dt\le \lambda C.
$$

Hence there exists a sequence of times $T_n\to\infty$ such
\begin{equation}
\label{eq2.144}\sup_n\left\|S_{\lambda}f(\cdot, T_n)\right\|_{L^2}^2<+\infty.
\end{equation}
The second and last step is to obtain a global space-time estimate for $f$ and therefore for $u$. For proving it  we use again identities \eqref{eq2.4} and \eqref{eq2.5} with $F= a(u) (f_x-\lambda\varphi^{\prime}(x-bt)f)$ and with $\eta$ defined on $[0,T_n]$ regular and positive, and  such that $\eta(0)=0$ and $\eta\equiv1$ if $t\ge 1$. With this choice of $\eta$ the right hand side of \eqref{eq2.4} is bounded using \eqref{eq1.5} and \eqref{eq2.12}. For the right hand side of \eqref{eq2.5} we use \eqref{eq2.14} and \eqref{eq2.144}. Hence we get that 
$$\frac{A_0}{2}\D\int_1^{\infty}\int
\varphi^{\prime\prime}_0\left(\lambda^3f^2+\lambda^2f^2_x\right)\,dxdt<+\infty,
$$
and as a consequence
$$\D\int_1^{\infty}\int e^{2\lambda|x-bt|}|u(x,t)|^2 \,dx dt< +\infty.$$
This implies that there exists a sequence of times $t_n$ such that
$$\lim_{t_n\to\infty}\|u(\cdot,t_n)\|=0.$$
But from the $L^2$ conservation law $u_0\equiv 0$. The result follows from the uniqueness property of the initial value problem in $H^1$ given in \cite{KPV1}.

In order to complete the proof of Theorem \ref{Th1} it remains to prove Claim \ref{CL1}.

\underline{Proof of Claim \ref{CL1}} We shall consider five different cases depending on the sign of $b$ and the relation between $b$ and $\lambda^2$. In what follows $C_0$ will denote a universal constant that can change from line to line. Only a finite number of choices of these constants will be made.

\underline{\sc Case \rm 1}. $\D\frac b{3\lambda^2}\ge 1/2$. In this case we shall take $\varphi=\D\frac 14\varphi_0$.

Recall \eqref{eq2.1} and \eqref{eq2.2}. Then, we have $\left(\varphi^{\prime}_0\le2\right)$

\begin{equation}
\label{eq2.15}
\textcircled{\tiny{5}}+\textcircled{\tiny{6}}+\textcircled{\tiny{12}}\ge\D\frac 1{16}\lambda b^2\int\varphi^{\prime\prime}f^2.
\end{equation}

On the other hand
\begin{equation}
\label{eq2.16}
\textcircled{\tiny{7}}+\textcircled{\tiny{8}}+\textcircled{\tiny{9}}+\textcircled{\tiny{10}}+\textcircled{\tiny{11}}
\le\D\int\varphi^{\prime\prime}f^2
\left\{C_0\lambda^3+\lambda C_0+\lambda bC_0\right\}.
\end{equation}

Hence taking $\lambda^{I}_0$ large enough we get from \eqref{eq2.15}--\eqref{eq2.16} and for $\lambda>\lambda^I_0$ a bound for
\begin{equation}
\label{eq2.17}
\lambda b^2\D\int\varphi^{\prime\prime}f^2=\frac{\lambda b^2}4\int\varphi^{\prime\prime}_0f^2\ge
\frac 9{16}\lambda^5\int\varphi^{\prime\prime}_0f^2.
\end{equation}

We also have
$$\textcircled{\tiny{4}}\le 6C_0\lambda\D\int\varphi^{\prime\prime}f_x^2.$$

Then using $\textcircled{\tiny{2}}$ and taking $\lambda^I_0$ large enough we get a bound of
\begin{equation}
\label{eq2.18}
\lambda b\D\int\varphi^{\prime\prime}f_x^2.
\end{equation}

Finally $\textcircled{\tiny{1}}$ gives a bound for $\left(f_{xx}\right)^2$ and $\textcircled{\tiny{3}}$ is positive but worse than $\textcircled{\tiny{2}}$. Putting everything together we get for $\lambda>\lambda_0^I$ that
\begin{equation}
\label{eq2.19}
\begin{array}{rcl}
\langle\left[S_{\lambda};A_{\lambda}\right]f,f\rangle&\ge&
A_0^I\D\int\varphi^{\prime\prime}_0\left(\lambda(f_{xx})^2+\lambda bf_x^2+\lambda b^2f^2\right) dx\\\\
&\ge&A_0^I\D\int\varphi_0^{\prime\prime}\left(\lambda^3f^2+\lambda^2f_x^2\right)dx,
\end{array}
\end{equation}
for some $A_0^I>0.$

\underline{\sc Case \rm 2}. $0\le\D\frac b{3\lambda^2}\le \frac12$, and the integrals in $\textcircled{\tiny{1}}, \textcircled{\tiny{2}},...\textcircled{\tiny{12}}$ will be considered in the region $r\ge 1$. We take $\varphi=\varphi_0$.

Then
\begin{equation}
\label{eq2.20}
\textcircled{\tiny{5}}+\textcircled{\tiny{6}}+\textcircled{\tiny{12}}\ge\D\frac 94\lambda^5\int_{r\geq1}\varphi^{\prime\prime}_0f^2;
\end{equation}
while the absolute value of the integral in the region $r\ge 1$ of  $\textcircled{\tiny{7}}+\textcircled{\tiny{8}}+\textcircled{\tiny{9}}+\textcircled{\tiny{10}}+\textcircled{\tiny{11}}$ is upper bounded by
\begin{equation}
\label{eq2.21}
\D\int_{r\ge 1}\varphi^{\prime\prime}_0f^2\left(C_0\lambda^3+C_0\lambda+C_0\lambda b\right).
\end{equation}

Hence taking $\lambda_0\ge\lambda_0^{II}$ for some $\lambda_0^{II}$ late enough we get from \eqref{eq2.20}--\eqref{eq2.21} a bound for
\begin{equation}
\label{eq2.22}
\lambda^5\D\int_{r\ge 1}\varphi^{\prime\prime}_0f^2.
\end{equation}

Using this time $\textcircled{\tiny{3}}$ instead of $\textcircled{\tiny{2}}$ (that is also positive) and $\textcircled{\tiny{1}}$ and \eqref{eq2.22} we get for the region $r\ge1$ a lower bound for the commutator of the type
\begin{equation}
\label{eq2.23}
A_0^{II}\D\int_{r\ge 1}\varphi^{\prime\prime}_0\left(\lambda(f_{xx})^2+\lambda^3f_x^2+\lambda^5f^2\right) dx,
\end{equation}
that is better than what we need.

\underline{\sc Case \rm 3}. $0\le\D\frac b{3\lambda^2}\le\frac 12$ and  $r\le 1$. We take $\varphi=\varphi_0$.

We need a cut off function $\eta$ such that $\eta(r)=r$ if $r\le 1$, $\eta(r)\equiv 0$ if $r\ge 3/2$ and $\left|\D\frac{d^k\eta}{dr^k}\right|\le C_0$ if $k=0,1,2$. Then for $r=|x-bt|$
$$\begin{array}{rcl}
\D\int\eta(r) f^2(x) dx&=&\D\frac 1{6\lambda^2}\int\left(3\lambda r^2-b\right)\left(\eta f^2\right)_{xx} dx\\\\
&=&\D\frac 1{6\lambda^2}\int\left(3\lambda r^2-b\right)
\left\{\left(2ff_{xx}+2f_x^2\right)\eta+\eta^{\prime\prime}f^2+2\eta^{\prime}ff^{\prime}\right\}.
\end{array}$$

Notice that $b\ge 0$, so we get, using Cauchy--Schwartz,
\begin{equation}
\label{eq2.24}
\begin{array}{rl}
\D\int\eta\left(f^2+\frac b{3\lambda^2}f_x^2\right)&\\\\
&\le\D\frac 1{6\lambda^2}\Bigg(\frac 12\int\left(3\lambda^2 r^2-b\right)^2 f^2\eta+2\int\left(f_{xx}\right)^2\eta+6\lambda\int r^2f_x^2\eta\\\\
&+C_0\D\int_{1\le r\le 3/2} f^2+\lambda f^2+\frac 1{\lambda} f_x^2\Bigg).
\end{array}
\end{equation}

From the definition of $\varphi_0$ we have that the terms $\textcircled{\tiny{4}}$, $\textcircled{\tiny{8}}$, $\textcircled{\tiny{9}}$, $\textcircled{\tiny{10}}$ and $\textcircled{\tiny{11}}$ are zero and that for $r\le 3/2$ we have that
$\textcircled{\tiny{5}}+\textcircled{\tiny{6}}+\textcircled{\tiny{12}}$ gives
$$\lambda\D\int_{r\le 3/2}\left(3\lambda^2r^2-b\right)^2\varphi^{\prime\prime}_0f^2.$$
Therefore we just have to take care of $\textcircled{\tiny{7}}$. Using \eqref{eq2.24} we get
$$\begin{array}{rcl}
3\lambda^3 \D\int_{r\le 1} f^2&\le& 3\lambda^3\D\int f^2\eta\\\\
&\le&\D\frac 14\lambda\int_{r\le 3/2}\left(3\lambda r^2-b\right)^2 f^2+
\lambda\int_{r\le 3/2}\left(f_{xx}\right)^2+3\lambda^2\int r^2 f_x^2 \eta\\\\
&&+C_0\left(\lambda^4\D\int_{1\le r\le 3/2} f^2+\lambda^2\int_{1\le r\le 3/2} f_x^2\right)\\\\
&=&I_1+I_2+I_3+I_4.
\end{array}$$

$I_1$ is bounded by $\textcircled{\tiny{5}}+\textcircled{\tiny{6}}+\textcircled{\tiny{12}}$ because $\varphi_0^{\prime\prime}=1$ if $r\leq\frac32$, $I_2$ by $\textcircled{\tiny{1}}$ , $I_3$ by $\textcircled{\tiny{3}}$, and $I_4$ by the bound \eqref{eq2.23} obtained in Case 2, and that as we said was better than what we needed in terms of powers of $\lambda$. As a consequence for $\lambda\ge \lambda_0^{III}$ and for some $\lambda_0^{III}$ large enough we get a lower bound for the commutator of the type
\begin{equation}
\label{eq2.25}
\D\int_{r\le 1}\varphi^{\prime\prime}_0\left(\lambda^3 f^2+\lambda(f_{xx})^2\right) dx.
\end{equation}

Finally from
$$\lambda^2\D\int f_x^2 \eta=-\lambda^2\int (f_x f \eta^{\prime}+ f f^{\prime\prime} \eta)=
\lambda^2\int(\frac 12 f^2 \eta^{\prime\prime}-f f^{\prime\prime} \eta),$$
and \eqref{eq2.25} we prove that there exists $A_0^{III}$ and a lower bound of the type
\begin{equation}
\label{eq2.255}
A_0^{III}\D\int_{r\le 1}\varphi_0^{\prime\prime}\left(\lambda^3 f^2+\lambda^2 f_x^2+\lambda(f_{xx})^2\right).
\end{equation}

\underline{\sc Case \rm 4}. $b\le 0$ and $ r\ge 1$. We take $\varphi=\varphi_0$.

If $b\le 0$, $\textcircled{\tiny{2}}$ comes with the wrong sign while $\textcircled{\tiny{12}}$ appears with the good one. Therefore we need another identity.
We have
\begin{equation}
\label{eq2.26}
\begin{array}{rcl}
\textcircled{\tiny{1}}+\textcircled{\tiny{2}}+\textcircled{\tiny{6}}&=&
9\lambda\D\int\varphi^{\prime\prime}_0(f_{xx})^2-
6\lambda b\left(\int\varphi_0^{\prime\prime} f f_{xx} +\varphi_0^{\prime\prime\prime} f_x f\right)+
\lambda b^2\int\varphi_0^{\prime\prime} f^2\\\\
&=&9\lambda\D\int\varphi^{\prime\prime}_0\left(f_{xx}-\frac b3 f\right)^2+3\lambda b\int \varphi^{IV}_0 f^2.
\end{array}
\end{equation}

Notice that for $r\ge 1$ the absolute value of the last term of the above identity
$$|\,3\lambda b\D\int \varphi^{IV}_0 f^2\,|,$$
together with the absolute values of $\textcircled{\tiny{7}}$, $\textcircled{\tiny{8}}$, $\textcircled{\tiny{9}}$, $\textcircled{\tiny{10}}$, and $\textcircled{\tiny{11}}$ are upper bounded by
$$C_0\D\int_{r\ge 1}\varphi_0^{\prime\prime} f^2\left(\lambda|b|+\lambda^3+\lambda\right).$$

But from $\textcircled{\tiny{5}}$ and $\textcircled{\tiny{12}}$ we get a lower bound of the type
$$\left(\lambda^5+|b|\lambda^3\right)\D\int_{1\leq r }\varphi_0^{\prime\prime} f^2,$$
that  can be taken in the region $1\leq r\leq\frac32$ as 
$$\left(9\lambda^5+6|b|\lambda^3\right)\D\int_{1\leq r\leq\frac32 }\varphi_0^{\prime\prime} f^2.$$

Also $\textcircled{\tiny{3}}$ gives the bound needed for $\textcircled{\tiny{4}}$. Altogether we get for $\lambda\geq\lambda^{IV}$  with $\lambda^{IV}$ large enough, that the commutator is lower bounded by
\begin{equation}
\label{eq2.27}
A_0^{IV}\D\int_{r\ge \frac32}\varphi_0^{\prime\prime}\left(\left(\lambda^5+|b|\lambda^3\right) f^2+\lambda^3 f_x^2\right)
+\D\int_{1\leq r\leq\frac32}\varphi_0^{\prime\prime}\left(\left(8\lambda^5+5|b|\lambda^3\right) f^2+17\lambda^3 f_x^2\right)
\end{equation}
for some $A_0^{IV}$ and that is better than what we need. 

\underline{\sc Case \rm 5}. $b\le 0$ and $r\le 1$. We take $\varphi=\varphi_0$.

As in Case 3 we need the cut off $\eta\equiv 1$ if $r\le 1$ and $\eta\equiv 0$ if $r\ge 3/2$. We have $r=|x-bt|$
$$\begin{array}{rcl}
\D\int f_x^2 \eta&=&
-\D\int r \partial_x\left(f_x^2\right) \eta-\int r f_x^2 \eta^{\prime}\\\\
&=&-\D\int 2r f_x f_{xx} \eta-\int r f_x^2 \eta^{\prime}\\\\
&=&-2\D\int\eta r\left(f_{xx}-\frac b3 f\right) f_x-
2\frac b3\int \eta rf f_x-\int r \eta^{\prime}f_x^2.
\end{array}$$

Hence using Cauchy--Schwarz and some integration by parts
\begin{equation}
\label{eq2.28}
\D\int f_x^2 \eta-\frac b3\int f^2 \eta
\le\D\frac 1{\lambda}\int\eta\left(f_{xx}-\frac b3 f\right)^2+\lambda\int\eta r^2 f_x^2+
\int r|\eta^{\prime}|\left(\frac{|b|}3 f^2+f_x^2\right).
\end{equation}

Also for $r\le 1$ the terms $\textcircled{\tiny{4}}$, $\textcircled{\tiny{8}}$, $\textcircled{\tiny{9}}$, $\textcircled{\tiny{10}}$, $\textcircled{\tiny{11}}$ are zero. Therefore only $\textcircled{\tiny{7}}$ remains as a negative term. We have
$$\begin{array}{rcl}
\D\int f^2 \eta&=&
\D\frac 12\int r^2(f^2)_{xx} \eta-2\int r \eta^{\prime} f f_x-\frac 12\int r^2 \eta^{\prime\prime} f^2\\\\
&=&\D\int r^2\left(f_x^2+f f_{xx}\right)\eta-2\int r\eta^{\prime} f f_x-\frac 12\int r^2 \eta^{\prime\prime} f^2\\\\
&=&\D\int r^2 f_x^2 \eta+\int r^2 f\left(f_{xx}-\frac b3 f\right) \eta+\frac b3\int r^2 f^2 \eta\\\\
&&\hphantom{\D\int r^2 f_x^2 \eta+\int r^2 f}-2\D\int r \eta^{\prime} f f_x-\frac 12\int r^2 \eta^{\prime\prime} f^2.
\end{array}$$

Then, using Cauchy--Schwarz
$$\begin{array}{rcl}
\D\int f^2\eta - \frac b3\int r^2 f^2 \eta
&\le&\D\int r^2\eta f_x^2+\frac{\lambda^2}2\int r^4 \eta f^2+\frac 1{2\lambda^2}\int\eta\left(f_{xx}-\frac b3 f\right)^2\\\\
&&+C_0\int_{1\le|r|\le 3/2}(1+\lambda)f^2+\frac 1{\lambda} f_x^2.
\end{array}$$

Then, for $r\le 1$ the absolute value of $\textcircled{\tiny{7}}$ is upper bounded by
$$\begin{array}{rcl}
3\lambda^3\D\int_{r\le 1}f^2\eta - \lambda^3 b\int r^2 f^2&\le&
3\lambda^3\D\int r^2 \eta f_x^2+\frac 32\lambda^5\int r^4 \eta f^2+\frac 32\lambda\int \eta\left(f_{xx}-\frac b3 f\right)^2\\\\
&&C_0\int_{1\le|r|\le3/2}\lambda^4f^2+\lambda^2f_x^2\\\\
&=&I_1+I_2+I_3+I_4.
\end{array}$$

Hence, we bound $I_1$ with $\textcircled{\tiny{3}}$. Notice that even in the region $1\le|r|\le 3/2$ there is no problem using \eqref{eq2.27}. For $I_2$ we use $\textcircled{\tiny{5}}$ and \eqref{eq2.27} because $\D\frac 32 r^4<8$ if $r\le 3/2$. For $I_3$ we use \eqref{eq2.26} $\left(\varphi_0^{IV}=0\,\text{ if }\, r\le 3/2\right)$. Finally for $I_4$ we use \eqref{eq2.27}.

As a conclusion we get a lower bound for the commutator of the type
$$\lambda^3\D\int_{r\le 1}f^2+ \lambda^3 |b|\int r^2 f^2+\lambda\int\eta\left(f_{xx}-\frac b3 f\right)^2+\lambda^3\int\eta r^2 f_x^2.$$

Hence using \eqref{eq2.28} we get a bound of $\lambda^2\D\int f_x^2 \eta$. Therefore there exists $A_0^V$ such that the commutator is bounded below by
\begin{equation}
\label{eq2.29}
A_0^V\D\int_{r\le1}\varphi_0^{\prime\prime}\left(\lambda^3 f^2+\lambda^2 f_x^2\right).
\end{equation}

Taking $A_0=\max\left\{A_0^I, A_0^{II}, A_0^{III}, A_0^{IV}, A_0^V\right\}$ we conclude the proof of Claim \ref{CL1} using \eqref{eq2.19}, \eqref{eq2.23}, \eqref{eq2.25}, \eqref{eq2.27} and \eqref{eq2.29}. Notice that a bound of $\lambda\D\int \varphi_0^{\prime\prime}(f_{xx})^2$ is also obtained.

\end{document}